\documentclass{amsart}

\newtheorem{theorem}{Theorem}[section]
\newtheorem{definition}[theorem]{Definition}
\newtheorem{proposition}[theorem]{Proposition}
\newtheorem{lemma}[theorem]{Lemma}

\newtheorem{corollary}[theorem]{Corollary}

\newcommand{\ZZ}{{\mathbb Z}}
\newcommand{\RR}{{\mathbb R}}
\newcommand{\CC}{{\mathbb C}}
\newcommand{\TT}{{\mathbb T}}

\newcommand{\DD}{{\mathbb D}}
\newcommand{\PP}{{\mathbb P}}
\newcommand{\QQ}{{\mathbb Q}}

\newcommand{\tr}{{\rm{tr}} }

\newcommand{\lsup}{{{\limsup}_{x\to x_0}}}

\title[On generalized sum rules for Jacobi matrices]
{On generalized sum rules for Jacobi matrices}

\author[F. Nazarov, F. Peherstorfer, A. Volberg and P. Yuditskii]
{F. Nazarov, F. Peherstorfer$^1$, A. Volberg and P. Yuditskii$^2$}
\thanks{$^{1,2}$ The work was supported by the Austrian Science
Found FWF, project number: P16390--N04}


\date{December 3, 2003}

\address{Department of Mathematics,
Michigan State University,
East Lansing, MI 48824}
\email{fedja@math.msu.edu}

\address{Institute for Analysis,
Johannes Kepler University of Linz, 
A-4040 Linz, Austria}
\email{Franz.Pehersorfer@jku.at}

\address{Department of Mathematics,
Michigan State University,
East Lansing, MI 48824}
\email{volberg@math.msu.edu}

\address{Institute for Analysis,
Johannes Kepler University of Linz, 
A-4040 Linz, Austria}
\email{Petro.Yudytskiy@jku.at}

\begin{document}

\begin{abstract}
		 This work is in a stream
		 (see e.g. \cite {DKS}, \cite {HS},
		 \cite {Ku}, \cite {LNS}, \cite{DeKu})
		 initiated by a paper of Killip and Simon \cite {KS},
   an earlier paper \cite{DK} also should be mentioned here. 
   Using methods of
Functional Analysis and the classical Szeg\"o Theorem 
we prove
sum rule identities in a very general form. Then, we apply 
the result to
obtain new asymptotics for orthonormal polynomials.
   \end{abstract}

\maketitle

\section{Introduction}
\subsection{Finite
dimensional perturbation 
of the 
Chebyshev matrix} 
Let $\{e_n\}_{n\ge 0}$ be the standard basis in
$l^2(\ZZ_+)$. Let $J$ be a Jacobi matrix
defining a bounded self--adjoint operator
on $l^2(\ZZ_+)$:
\begin{equation*}
Je_n=p_ne_{n-1}+q_n e_n +p_{n+1} e_{n+1},
\quad n\ge 1,
   \end{equation*}
   and
   $$
   Je_0=q_0 e_0 +p_{1} e_{1}.
   $$
 Under the condition  $p_n>0$, the vector $e_0$ is cyclic
 for $J$. The function 
 $$
 r(z)=\langle(J-z)^{-1}e_0,e_0\rangle
 $$
 is called the resolvent function. It has the representation
 $$
 r(z)=\int\frac{d\sigma(x)}{x-z}.
 $$
 The measure $\sigma$, $d\sigma\ge 0$, is called the spectral
 measure of $J$. 
 
 Using a three term recurrence relation
 for orthonormal polynomials $\{P_n(z)\}_{n\ge 0}$
 with respect to $\sigma$ one can restore the coefficient
 sequences of $J$ 
 \begin{equation*}
zP_n(z)=p_n P_{n-1}(z)+q_n P_n(z) +p_{n+1} P_{n+1}(z),
\quad n\ge 1,
   \end{equation*}
   and
   $$
   zP_0(z)=q_0 P_0(z) +p_{1} P_{1}(z).
   $$
   
  With a given $J$ we associate a sequence $J(n)$
defined
   by
   \begin{equation*}
p(n)_k=\begin{cases} p_k, &\ k<n
\\ 1, &\ k\ge n
\end{cases},
\end{equation*}
   \begin{equation*}
q(n)_k=\begin{cases} q_k, &\ k<n
\\ 0, &\ k\ge n
\end{cases}.
\end{equation*}
   $J(n)$ is a finite
dimensional perturbation of the ``free"
(Chebyshev) matrix 
$J_0=S_++S_+^{*}$, $S_+e_n=e_{n+1}$.

Note that
$$
r_0(z)=\langle(J_0-z)^{-1}e_0,e_0\rangle=
-\zeta,
$$
where $1/\zeta+\zeta=z,\ \zeta\in\DD$, that is
$\zeta=\frac{z-\sqrt{z^2-4}}{2}$.
Further, in terms of orthonormal polynomials
$$
r(n)(z)=\langle(J(n)-z)^{-1}e_0,e_0\rangle=-
\frac{p_n Q_n(z)-\zeta Q_{n-1}(z)}{p_n P_n(z)-\zeta P_{n-1}(z)},
$$
where $Q_n$ are so called orthonormal polynomials
of the second kind
$$
Q_n(z):=\int\frac{P_n(x)-P_n(z)}{x-z}d\sigma.
$$
(They satisfy the same three term recurrence relation 
as $P_n$'s but with a
different initial condition). What is important for us
\begin{equation}\label{add0}
\sigma'(n)_{a.c.}(x)=
\frac 1 \pi{\rm Im}r(n)(x+i0)=
\frac 1 \pi\frac{-{\rm Im}\zeta(x+i0)}
{|p_n P_n(x)-\zeta(x+i0) P_{n-1}(x)|^2}.
\end{equation}
and
\begin{equation}\label{add1}
\sigma(J(n))\cap\{\RR\setminus[-2,2]\}=
\{z\in\CC\setminus[-2,2]:p_n P_n(z)-\zeta(z) P_{n-1}(z)=0\}.
\end{equation}

The perturbation determinant of $J(n)$ with respect
to $J_0$ is well defined and
we can introduce a function 
\begin{equation*}
\Delta_n(\zeta)=\frac 1 {\prod_{j=1}^{n-1} p_j}\det(J(n)-z)
(J_0-z)^{-1}.
\end{equation*}
By definition
\begin{equation}\label{1}
\log \Delta_n(z)=
- t(n)_0-\sum_{k\ge 1}
 \frac{t(n)_k}{k z^k} 
\end{equation}
where
\begin{equation*}
t(n)_0=\sum_{j=1}^{n-1}\log p_j,\quad
t(n)_k=\tr (J(n)^k-J_0^k),\ 
k\ge 1.
\end{equation*}
On the other hand one can find the determinant 
by a direct calculation and then
$$
\Delta_n(z)=(p_n P_n(z)-\zeta P_{n-1}(z))\zeta^n,
$$
as before $1/\zeta+\zeta=z,\ \zeta\in\DD$.

Therefore $\Delta_n(z)$ 
has explicit representation \eqref{1} in terms of coefficients of
$J(n)$, on the other hand
it has nice analytic properties:
   its zeros 
in $\overline\CC\setminus[-2,2]$ are simple  
and related to the eigenvalues of $J(n)$ in this region
(see \eqref{add1});
it has no poles; and  by \eqref{add0}
\begin{equation}\label{2}
|\Delta_n (x+i0)|^2=\frac 1{2\pi}
\frac{\sqrt{4-x^2}}{\sigma'(n)_{a.c.}}.
\end{equation}
That is, we can restore $\Delta_n(z)$ only in terms of these (partial)
spectral data (see the next subsection).

\subsection{The Killip--Simon functional via spectral data}
\begin{definition} Let $J$ be a Jacobi matrix with a spectrum on
$[-2,2]\cup X$, where the only possible accumulation points of
$X=\{x_k\}$ are
$\pm 2$. Following to Killip and Simon, to a given nonnegative
polynomial $A$ we associate the functional
 that might diverge only to
$+\infty$
\begin{equation}\label{3}
\Lambda_A(J):=
\sum_
{X}
F(x_k)+\frac 1 {2\pi}\int_{-2}^2
\log\left(
\frac{\sqrt{4-x^2}}{2\pi\sigma'_{a.c.}}\right)
A(x)\sqrt{4-x^2}\,dx,
\end{equation}
where
\begin{equation}\label{33.3}
\begin{split}
F(x)=&\int_2^x A(x)\sqrt{x^2-4}\,dx\quad
\text{for}\quad x>2,\\
F(x)=&-\int_{-2}^x A(x)\sqrt{x^2-4}\,dx\quad
\text{for}\quad x<-2.
\end{split}
\end{equation}
\end{definition}

Let us point out that 
the  Killip--Simon functional
$\Lambda_A(J)$ is defined in terms of the spectral data
of $J$ only.
Let us demonstrate how to obtain
for a finite dimensional perturbation $J(n)$ of $J_0$
a representation of $\Lambda_A(J(n))$ in terms of the recurrence
coefficients.

First, let us note that the function
$\log{\Delta_n(z)}$
is well defined in the upper half plane, in fact, in the domain
$\overline\CC\setminus \sigma(J(n))$. Moreover, the boundary
values of the real part ${\rm{Re}} \log\Delta_n(x+i0)$, $x\in [-2,2]$,
are given by \eqref{2}. For $x\ge 2$ the
imaginary part of $\log\Delta_n(z)$ (that is the
argument of $\Delta_n(\zeta)$) is of the form
$$
\frac 1\pi\arg \Delta_n(x+i0)=
\#\{y\in\sigma(J(n)):y\ge x\}
$$
and similarly, 
$$
\frac 1\pi\arg \Delta_n(x+i0)=-
\#\{y\in\sigma(J(n)):y\le x\}
$$
for $x\le -2$. 
Therefore, multiplying $\log\Delta_n(z)$ by
$A(z)\sqrt{z^2-4}$,  where $A(z)$ is 
the given nonnegative polynomial,
we get a function with 
the following representation
\begin{equation}\label{4}
A(z)\sqrt{z^2-4}\log\Delta_n(z)=
B_n(z)+\int_{\sigma(J(n))}\frac{d\lambda_n}{x-z},
\end{equation}
where $B_n(z)$ is a (real) polynomial of degree
one
bigger than $A$  and
\begin{equation*}
\lambda'_n(x)=
\begin{cases} \frac 1{2\pi}
A(x)\sqrt{4-x^2}\log 
\frac 1{2\pi}
\frac{\sqrt{4-x^2}}{\sigma'(n)_{a.c.}},
&x\in[-2,2]
\\ A(x)\sqrt{x^2-4}
\#\{y\in\sigma(J(n)):y\ge x\},
&x\ge 2\\
A(x)\sqrt{x^2-4}
\#\{y\in\sigma(J(n)):y\le x\},
&x\le -2
\end{cases}.
\end{equation*}
Thus the functional $\Lambda_A(J(n))=\int d\lambda_n$.

Let us mention that the polynomial $B_n(z)$
is determined uniquely by \eqref{4} since
\begin{equation}\label{25.1}
\int_{\sigma(J(n))}\frac{d\lambda_n}{x-z}
=-\frac{\int d\lambda_n}{z}-...=
\underline O\left(\frac 1 z\right),\quad
z\to\infty.
\end{equation}

Let us define
$$
\Phi(z)=\text{Const}+a_1 z+\dots+ a_{m+2}
z^{m+2}
$$
by
$$
\Phi'(z)=zA(z)-
\frac 1\pi\int_{-2}^2
\frac{A(x)-A(z)}{x-z}\sqrt{4-x^2}\,dx.
$$
Note that
\begin{equation}\label{5}
A(z)\sqrt{z^2-4}=
\frac 1\pi\int_{-2}^2
\frac{A(x)}{x-z}\sqrt{4-x^2}\,dx
+\Phi'(z).
\end{equation}
Therefore, using \eqref{1},
\eqref{4}, \eqref{25.1} and \eqref{5} we get
\begin{equation}\label{6}
\begin{split}
\int d\lambda_n=&
   - a t(n)_0+
a_1 t(n)_1+
2a_2\frac{t(n)_2}{2}+\dots+
(m+2) a_{m+2}\frac{t(n)_{m+2}}{(m+2)}\\
=&-a
   t(n)_0
+\tr\{\Phi(J(n))-\Phi(J_0)\},
\end{split}
\end{equation}
where we put
$$a=
\frac 1\pi\int_{-2}^2
{A(x)}\sqrt{4-x^2}\,dx.
$$

Note, if
$A(z)=1$, that is $a=2$,
$\Phi(z)=\text{Const}+z^2/2$, 
then we are in the Killip--Simon case \cite{KS}:
\begin{equation*}
\int d\lambda_n=\frac{t(n)_2} 2-2t(n)_0
=-\frac 1 2+\sum_{k=1}^{\infty}(p(n)_k^2-1-\log p(n)_k^2)
+\frac 1 2\sum_{k=0}^{\infty} q(n)^2_k.
\end{equation*}
For a more general example see Appendix.

\subsection{The Killip--Simon functional via coefficient sequences}
For a bounded operator $G$ in $l^2(\ZZ_+)$ we denote
$G^{(k)}:=(S^*_+)^k G S_+^k$.
\begin{lemma} For all $k\ge 1$ and $n\ge l-1$
$$
(J^{(k)})^le_n=(J^l)^{(k)}e_n.
$$
\end{lemma}
\begin{proof}
Let us mention that the decomposition of the vector $J^le_{k+n}$
begins with the basic's vector $e_{k+n-l}$. Therefore the
orthoprojector 
$P_{k-1}$ onto the
subspace spanned by $\{e_0,...e_{k-1}\}$ annihilates this vector,
$P_{k-1}J^le_{k+n}=0$. Thus, by induction,
\begin{equation*}
\begin{split}
(J^{(k)})^{l+1}e_n=&J^{(k)}(J^{(k)})^le_n=
J^{(k)}(J^l)^{(k)}e_n=
(S^*_+)^kJS_+^k(S^*_+)^kJ^lS_+^ke_n
\\=&(S^*_+)^kJ(I-P_{k-1})J^le_{k+n}=(S^*_+)^kJ^{l+1}e_{k+n}=
(J^{l+1})^{(k)}e_n.
\end{split}
\end{equation*}
\end{proof}

For a bounded Jacobi matrix $J$
(and a polynomial
$A$) let us define a function of a finite number of variables
$$
h_A=h_A(J):=-a\log p_{m+2}+\langle
\{\Phi(J)-\Phi(J_0)\}e_{m+1},e_{m+1}\rangle.
$$
Note that due to the previous lemma
\begin{equation*}
\begin{split}
h_A\circ\tau^k&=-a\log p_{m+k+2}+\langle \{\Phi(J^{(k)})-\Phi(J_0)\}
e_{m+1},e_{m+1}\rangle\\
&=-a\log p_{m+k+2}+\langle \{\Phi(J)-\Phi(J_0)\}
e_{m+k+1},e_{m+k+1}\rangle,
\end{split}
\end{equation*}
where $\tau$ acts just as a shift of indexes. In this case the series
\begin{equation*}
\sum_{k\ge 0}h_A\circ\tau^k
\end{equation*}
may not converge, but the generic term is well define.
\begin{definition}
With a given Jacobi matrix $J$
and a polynomial $A$ of degree $m$
we associate the series
\begin{equation}\label{s1}
H_A(J):=
\sum_{k=0}^{m}(-a\log p_{k+1}+\langle \{\Phi(J)-\Phi(J_0)\}
e_{k},e_{k}\rangle)
+\sum_{k\ge 0}h_A\circ\tau^k.
\end{equation}
\end{definition}
Note that
 $H_A(J(n))$ is just a
finite sum, in fact $h\circ\tau^k$ vanishes starting with
a suitable $k$, moreover $H_A(J(n))=\Lambda_A(J(n))$.

\subsection{Results}
\begin{theorem}\label{1.1}
	Let $A$ be a nonnegative polynomial.
 The spectral measure $\sigma$ of a
 Jacobi matrix $J$ with a spectrum
 of the form $[-2,2]\cup X$, where
 $\pm 2$ are the only possible accumulation
 points of the discrete set $X$, satisfies
$\Lambda_A(J)<\infty$ if and only if  series
\eqref{s1} converges;
moreover $H_A(J)=\Lambda_A(J)$. 
\end{theorem}

In a sense our result is a kind of ``existence theorem".
 To balance the situation we derive from it the following
 application. (We conjectured this result in a note mentioned
 in \cite{Ku}).
 
 \begin{theorem}\label{A}
 	Let $A(x)$ be a nonnegative polynomial
  of degree $m$
  with all zeros on $[-2,2]$.
 Let a measure $\sigma$ supported on $[-2,2]\cup X$
 satisfy the condition $\int d\lambda<\infty$, where
 \begin{equation}\label{lambda}
\lambda'(x)=\lambda'(x;\sigma)=
\begin{cases} \frac 1{2\pi}
A(x)\sqrt{4-x^2}\log\left( 
\frac 1{2\pi}
\frac{\sqrt{4-x^2}}{\sigma'_{a.c.}(x)}\right),
&x\in[-2,2]
\\ A(x)\sqrt{x^2-4}
\#\{y\in X:y\ge x\},
&x\ge 2\\
A(x)\sqrt{x^2-4}
\#\{y\in X:y\le x\},
&x\le -2
\end{cases}.
\end{equation}
 Then
 the sequence of
  orthonormal polynomials
 $P_n(z)=P_n(z;\sigma)$, 
 normalized by
 \begin{equation*}
\zeta^{n+1}\sqrt{z^2-4} P_n(z)\exp\left(-\frac{\tilde
B_n(z)}{A(z)\sqrt{z^2-4}}
\right)=1+\underline O\left(\frac 1{z^{m+2}}\right),
 \end{equation*}
the polynomial $\tilde B_n(z)$ (of degree ${m+1}$)
 is determined uniquely by the condition
  \begin{equation*}
\log\{\zeta^{n+1}\sqrt{z^2-4} P_n(z)\}-\frac{\tilde
B_n(z)}{A(z)\sqrt{z^2-4}} =\underline O\left(\frac 1{z^{m+2}}\right),
 \end{equation*}
 converges uniformly on compact subsets of the 
 domain $\overline\CC\setminus [-2,2]$ to the holomorphic function
 \begin{equation}\label{D}
D(z):= 
\exp\left(
\frac{1}{A(z)\sqrt{z^2-4}}
\int\frac{d\lambda}{x-z}
\right).
\end{equation}
\end{theorem}
Note that
as well as in the Szeg\"o case
the limit function $D(z)$ can be expressed
only in terms of
 $\sigma'_{a.c.}$ and $X$.

\section{Semicontinuity of Szeg\"o type functional}

For a measure $\mu$ on the unit circle $\TT$ we denote
by ${\rm Sz}(\mu)$ the functional
$$
{\rm Sz}(\mu)=\int_{\TT}\log\frac{d\mu_{a.c.}}{dm}\,dm.
$$
Recall the main property of this functional
 $$
{\rm Sz}(\mu)=
\inf\{\log\int_\TT |1-f|^2\,d\mu(t):
f \ \text{is a polynomial}, f(0)=0\}.
$$

\begin{lemma}\label{ls1}
Let $\mu_k$ converge weakly to  $\mu$. Then
\begin{equation}\label{sz1}
\limsup {\rm Sz}(\mu_k)\le {\rm Sz}(\mu).
\end{equation}
\end{lemma}

\begin{proof}
Since
for every $\epsilon$ there exists a polynomial
$g$, $g(0)=0$, such that
$$
\log\int_\TT |1-g|^2\,d\mu(t)\le {\rm Sz}(\mu)+\epsilon,
$$
starting from a suitable $k$ we have
$$
\log\int_\TT |1-g|^2\,d\mu_k(t)\le
{\rm Sz}(\mu)+2\epsilon.
$$
But for every $k$
\begin{equation*}
\begin{split}
{\rm Sz}(\mu_k)=&
\inf\{\log\int_\TT |1-f|^2\,d\mu_k(t):
f\ \text{is a polynomial}, f(0)=0\}\\
\le&
\log\int_\TT |1-g|^2\,d\mu_k(t).
\end{split}
\end{equation*}
Thus \eqref{sz1} is proved.
\end{proof}

\begin{lemma}\label{2.2}
Let $\rho$ be a normalized nonnegative weight, i.e., $\rho\ge 0$, 
$\int_{\TT}\rho\,dm=1$,  such that $\rho\log\rho\in L^1$.
Assume that $\mu_k$ converges weakly to  $\mu$. Then
\begin{equation}\label{s2}
\liminf \int_{\TT}\log\frac{dm}{d(\mu_k)_{a.c.}}\,\rho dm
\ge \int_{\TT}\log\frac{dm}{d\mu_{a.c.}}\,\rho dm.
\end{equation}
\end{lemma}

\begin{proof}
Define a map $\psi:\TT\to\TT$ by
$\psi(e^{i\theta})=\exp\{i\int_0^{\theta}\rho(e^{i\theta})\,d\theta\}$
and denote by $\phi$ the inverse map, 
 $\psi\circ\phi={\rm id}:\TT\to\TT$. Let us apply Lemma \ref{ls1}
 to the sequence $\tilde\mu_n:=\mu_n\circ\phi$ that converges
 weakly to $\tilde\mu:=\mu\circ\phi$.
 \begin{equation*}
\liminf \int_{\TT}\log\frac{dm}{d(\tilde\mu_k)_{a.c.}}\, dm
\ge \int_{\TT}\log\frac{dm}{d\tilde\mu_{a.c.}}\, dm.
\end{equation*}
Making the inverse change of variable in each integral
we have
 \begin{equation*}
\liminf \int_{\TT}\log\frac{\rho dm}{d(\mu_k)_{a.c.}}\, \rho dm
\ge \int_{\TT}\log\frac{\rho dm}{d\mu_{a.c.}}\,\rho dm.
\end{equation*}
Since  $\rho\log\rho\in L^1$ we get \eqref{s2}.

\end{proof}

\begin{corollary}\label{c2.3}
$$
\liminf_{n\to \infty}\Lambda_A(J(n))\ge\Lambda_A(J).
$$
\end{corollary}
\begin{proof}
Outside of $[-2,2]$ we apply the Fatou Lemma, e.g. 
\cite{Yo}, p. 17, and
on $[-2,2]$ we apply Lemma \ref{2.2}
\end{proof}

\section{ Lemma on positiveness and its consequences}

For a given interval $I$, $0\in I$, let
$h\in C(I^l)$ be such that
$h(0,...,0)=0$. Then
$$
H(\underline x)=\sum_{i=0}^\infty
h(x_{i+1},x_{i+2},...,x_{i+l})
$$
is well defined on
$$
I^\infty_0=
\{\underline x:\underline x=
(x_{0},x_{1},...,x_{n},0,0...)\}.
$$

\begin{lemma}\label{ML}
Assume that
$H$ is bounded from below, $H(\underline x)\ge C$ for all
$\underline x\in I^\infty_0$.
Then there exists a function $g$ of the form
$$
g(x_{1},...,x_{l})=
h(x_{1},...,x_{l})+
\gamma(x_{2},...,x_{l})-
\gamma(x_{1},...,x_{l-1}),\quad
\gamma\in C(I^{l-1}),
$$
such that $g\ge 0$.
\end{lemma}

First we prove a sublemma.

\begin{lemma} The set $G$, consisting of functions of
the form
\begin{equation*}
G=\{g(x_{1},...,x_{l})+\gamma(x_{1},...,x_{l-1})-
\gamma(x_{2},...,x_{l})\},
\end{equation*}
where
 $g\in C(I^l), g\ge 0, g(0)=0$,
$\gamma\in C(I^{l-1})$,
is closed in $C(I^l)$.
\end{lemma}

\begin{proof}

We give a proof in the case of two variables
 (the general case can be considered in a similar way).
 
 Let
 \begin{equation}\label{pos0}
 h(x,y)=\lim \{g_n(x,y) +\gamma_n(x)-\gamma_n(y)\},
 \end{equation}
Assuming the normalization $\gamma_n(0)=0$ we get a uniform bound
 for $\gamma_n$,
 $$
 -1-h(0,x)\le \gamma_n(x)\le h(x,0)+1.
 $$
 Therefore there exists a subsequence that converges
 weakly, say, in $L^2$. Then, using the Mazur Theorem, see 
 e.g. \cite{Yo}, p. 120, and convexity of $G$ we can find a sequence 
 $\gamma_n^{(1)}(x)$ and corresponding sequence of 
 $g_n^{(1)}(x,y)\ge 0$
 such that
 $\gamma_n^{(1)}\to \gamma_1$, $g_n^{(1)}\to g_1$ 
 in $L^2$ strongly and we still have \eqref{pos0}.

Thus, there exists a representation
 \begin{equation}\label{pos1}
h(x,y)=g_1(x,y)+\gamma_1(x)-\gamma_1(y)
 \end{equation}
that holds almost everywhere, and the function $\gamma_1(x)$,
in fact, because of uniform boundness, belongs to
$L^{\infty}$.

Starting with this place we will show that
there exists a representation for
$h(x,y)$ of the form \eqref {pos1} but with
continuous functions $\gamma$ and $g\ge 0$.
First, let us construct
a function
$\gamma_2$ which is defined for all $x\in I$ and such that
$\gamma_2(x)-\gamma_2(y)\le h(x,y)$ holds everywhere. 

Set 
$\gamma_2(x_0)=\limsup_{\delta\to 0}\frac 1
{2\delta}\int_{x_0-\delta}^{x_0+\delta}
\gamma_1$. Note that $\gamma_2(x)=\gamma_1(x)$ (a.e.). To show that
$\gamma_2(x)-\gamma_2(y)\le h(x,y)$ for all $(x,y)\in I^2$ we average
the inequality with $\gamma_1$ over rectangles 
$\{x_0-\delta\le x\le x_0+\delta,\ y_0-\delta\le y\le
y_0+\delta\}$ and take the upper limit when $\delta\to 0$. Since
$$
\limsup(a+b)\ge \limsup a+\liminf b
$$
we get the inequality we need.
Next, we construct an upper semicontinuous function
$\gamma_3(x_0)=\lsup\gamma_2(x)$.

Let $\Gamma$ be the set of upper semicontinuous functions
defined on $I$ with normalization $\gamma(0)=0$ and
such that $\gamma(x)-\gamma(y)\le h(x,y)$. The previous
construction shows that 
$\Gamma\not =\emptyset$. Now, the key point  is to consider the
function
$$
\gamma_4(x):=\sup\{\gamma(x): \gamma\in \Gamma\}.
$$
It belongs to $\Gamma$ since $\sup\{\beta_1(x),\beta_2(x)\}\in\Gamma$
if only $\beta_1(x)\in \Gamma, \beta_2(x)\in\Gamma$.

We claim that $\gamma_4(x)$ is {\it lower} semicontinuous. Assume, 
on the contrary, that it is not. This means that there exist $\delta
>0$,  a point $x_0\in I$ and a sequence $\{x_n\}$,
$\lim x_n=x_0$, such that $\gamma_4(x_n)\le\gamma_4(x_0)-\delta$.
Let us mention that $x_0\not= 0$ since
$$
-h(0,x)\le \gamma(x)\le h(x,0),
$$
and hence $\lim_{x\to 0}\gamma(x)=0=\gamma(0)$ for all
$\gamma\in\Gamma$.

The function $h(x,y)$ is continuous therefore we can choose
such $N$ that 
$$|h(x_N,y)-h(x_0,y)|\le \delta/2
$$ 
for all $y\in I$.

Let 
\begin{equation*}
\gamma_5(x)=\begin{cases} \gamma_4(x), &\ x\not=x_N 
\\ \gamma_4(x_N)+\delta/2, &\ x=x_N.
\end{cases}
\end{equation*}
Let us check that $\gamma_5\in \Gamma$. 
It is upper semicontinuous, $\gamma_5(0)=0$. Further,
for $y\not =x_N$ we
have
\begin{equation*}
\begin{split}
\gamma_5(x_N)-\gamma_5(y)=&\gamma_4(x_N)+\delta/2-\gamma_4(y)\\
\le &
\gamma_4(x_0)-\delta/2-\gamma_4(y)\\
\le & h(x_0,y)-\delta/2\le h(x_N,y).
\end{split}
\end{equation*}
Moreover the inequality $\gamma_5(x)-\gamma_5(y)\le h(x,y)$ holds on the
line
$y=x_N$ and for all other values of $x$ and $y$.

On the other hand it could not be in the class, since
$$
\gamma_5(x_N)>\sup\{\gamma(x_N),\ \gamma\in \Gamma\}.
$$
Therefore we arrive to a contradiction. Thus $\gamma_4(x)$ is
simultaneously upper and lower semicontinuous, that is 
$\gamma_4(x)$ is a continuous function. The lemma is proved.

\end{proof}

\begin{proof}[Proof of Lemma \ref{ML}]
If not then $h$ does not belong to the closed convex set $G$.
Therefore there exists a measure
$\mu\in C(I^l)^*$, $d\mu\ge 0$, such that 
\begin{equation}\label{l8}
\int_{I^l}h(x)\,d\mu(x)<0
\end{equation}
and
\begin{equation*}
\int_{ I^l}
(\gamma(x_{2},...,x_{l})-
\gamma(x_{1},...,x_{l-1}))
\,d\mu(x)=0.
\end{equation*}
In other words
\begin{equation}\label{l9}
\int_{z\in I}d\mu (y,z)=
\int_{z\in I}d\mu (z,y)
\end{equation}
for all $y\in I^{l-1}$.

Without lost of generality we may assume that $\mu$ 
is absolutely continuous,
moreover $d\mu =w(x_1,\ldots,x_l)\,dx_1\ldots dx_l$, $w\not =0$ a.e.
Note that condition \eqref{l9} is now of the form
\begin{equation}\label{9.1}
\int_{z\in I}w (y,z)\,dz=
\int_{z\in I}w(z,y)\,dz,
\quad y\in I^{l-1}.
\end{equation}

We want to get a contradiction between
\eqref{l8} and $H\ge C$ by
extending the functional related
to $w$ on functions on $I^\infty_0$.

We can normalize $w$ by the condition
$\int_{I^l} w(x)=1$. Let us think
on $w$ as on the probability
$$
w(y) dy=\PP\{\underline x: x_i\in (y_i,y_i+dy_i),\
i=1,...,l\},
$$
and we want
\begin{equation}\label{l10}
\PP\{\underline x: x_{i+k}\in (y_i,y_i+d y_i),\
i=1,...,l\}=w(y)\,dy,\
\text{for all}\ k,
\end{equation}
that is the probability should be shift
invariant.
Actually we will define
on $I^N$ step by
step for increasing $N$ 
probabilistic measures 
$$
\rho(x_1,\ldots,x_N) dx_1\ldots dx_N
$$
 using a conditional probability. 

For $N\ge l$ inductively define
\begin{equation*}
	\begin {split}
	&\rho(x_1,\ldots,x_N, x_{N+1}) dx_1\ldots dx_N dx_{N+1}:=\\&
	\rho(x_1,\ldots,x_N) dx_1\ldots dx_N
\frac{w(x_{N+2-l},...,x_N,x_{N+1})dx_{N+1}}
{\int_{ I} w(x_{N+2-l},...,x_N ,v)dv}.
	\end {split}
\end{equation*}
Now we have to check that \eqref{l10}
holds true.

If  $k\not=N+1-l$ then  \eqref{l10}  holds by the induction conjecture
since
$$
\int_{ I}\rho(x_1,\ldots,x_N, x_{N+1}) dx_{N+1}=
\rho(x_1,\ldots,x_N).
$$
In case $k=N+1-l$ we have
\begin{equation*}
	\begin {split}
	&\int\rho(x_1,\ldots, x_{N+1-l},y_1,\ldots y_l) dx_1\ldots  dx_{N+1-l}\\=&
\int\left(\int_{x\in I^{N-l}}	\rho(x, x_{N+1-l},y_1,\ldots y_{l-1}) dx
\right)
dx_{N-l+1}
\frac{w(y_{1},...,y_{l})}
{\int w(y_{1},...,y_{l-1},v)dv}\\=&
\int w(x_{N-l+1},y_{1},...,y_{l-1})
dx_{N-l+1}
\frac{w(y_{1},...,y_{l})}
{\int w(y_{1},...,y_{l-1},v)dv}.
	\end {split}
\end{equation*}
Making use of \eqref{9.1} we get 
$$
\int\rho(x_1,\ldots, x_{N+1-l},y_1,\ldots y_l) dx_1\ldots 
dx_{N+1-l}=
w(y_{1},...,y_{l})
$$
  that is \eqref{l10} is proved.

Now we are in a position to finish Lemma's proof.
For $\underline x$'s of the form $\underline x=
(x,0,...),\  x\in I^{N}$, we can integrate $H$
against $\rho$:
$$
\int_{x\in I^N}
H(\underline x)
\rho( x)\ge C.
$$
On the other hand using the definition of $H$ and the key property of
$\rho$
we get
\begin{equation}\label{l12}
C\le \int_{x\in I^N}
H(\underline x)
\rho( x)\le (l-1)||h||+(N-l+1)\int_{I^l} h(y)w(y)\,dy.
 \end{equation}
 Since $N$ is arbitrary large,
\eqref{l8} contradicts to \eqref{l12}.
\end{proof}

\begin{corollary}
For a nonnegative polynomial $A$ there exist
continuous
functions $g_A$ and $\gamma_A$ such that
\begin{equation}\label{c1}
h_A=g_A+\gamma_A\circ\tau-\gamma_A
\end{equation}
and $g_A\ge 0$.
\end{corollary}

\begin{proof}
Note that $H_A(J(n))$ are uniformly bounded from below.
\end{proof}

 \begin{corollary}
Let $J$ be such that $p_n\to 1$ and $q_n\to 0$. Then
\begin{equation}\label{c2}
H_A(J):=
\sum_{k=0}^{m}(-a\log p_{k+1}+\langle \{\Phi(J)-\Phi(J_0)\}
e_{k},e_{k}\rangle)-\gamma_A
+\sum_{k\ge 0}g_A\circ\tau^k.
\end{equation}
That is the series with positive terms $\sum_{k\ge 0}g_A\circ\tau^k$
converges if and only if the series $\sum_{k\ge
0}h_A\circ\tau^k$ converges. 
\end{corollary}
\begin{proof} We use representation \eqref{c1} and continuity of
$\gamma_A$.

 \end{proof}

 \section{Proof of the Main Theorem}
 
Assume that for a given $J$ its spectral measure
$\sigma$ is such that $\Lambda_A(J)<\infty$,
see definition \eqref{3}.
 Note that due to
Denisov--Rakhmanov Theorem \cite{D}
\begin{equation}\label{p1}
p_n(\sigma)\to 1,\quad
q_n(\sigma)\to 0
\end{equation}
and we can use \eqref{c2} as a definition of $H_A(J)$.

With the measure $\sigma$ let us associate a measure
$\sigma_\epsilon$ that we get by using 
the following two regularizations.
First, we add to its absolutely continuous part
the component $\epsilon\,dx$, that is 
$(\sigma'_\epsilon)_{a.c.}=\sigma'_{a.c.}+\epsilon$.
Second, we leave just a finite number of the spectral
points outside of $[-2,2]$, say, that one that belongs to
$\RR\setminus [-2-\epsilon, 2+\epsilon]$. It is important that
\begin{equation}\label{p2}
p_n(\sigma_\epsilon)\to p_n(\sigma),\quad
q_n(\sigma_\epsilon)\to q_n(\sigma)
\end{equation}
for a fixed $n$ as $\epsilon\to 0$.
The measure $\sigma_\epsilon$ satisfies the conditions of 
Szeg\"o's
Theorem, and therefore $\zeta^n P_n(z,\sigma_\epsilon)$ converges
uniformly on compact subsets of $\overline\CC\setminus[-2,2]$ to a
certain function that can be expressed directly in terms of
$(\sigma'_\epsilon)_{a.c.}$ and the mass--points 
outside of $[-2,2]$,
see e.g. \cite{PYu}. We use the consequence of this statement in the
form
\begin{equation*}
\zeta^n (p_n (\sigma_\epsilon) P_n(z,\sigma_\epsilon)-
\zeta P_{n-1}(z,\sigma_\epsilon))\to \Delta(z;\sigma_\epsilon)
\end{equation*}
uniformly on compact subsets of $\overline\CC\setminus[-2,2]$.
Here $\Delta(z;\sigma_\epsilon)$ is defined by
$$
\Delta(z;\sigma_\epsilon)=\exp\left\{
\sqrt{z^2-4}\int\frac{1}{x-z}
\frac{d\lambda(x;\sigma_\epsilon)}{x^2-4}\right\}.
$$
In other words
\begin{equation*}
\log\Delta(z;J(n;\sigma_\epsilon))
\to \log\Delta(z;\sigma_\epsilon),\ n\to \infty,
\end{equation*} 
 uniformly on $\overline\CC\setminus{\rm supp}(\sigma_\epsilon)$.
 
 Finally, since (all) coefficients in  decomposition
 \eqref{1} of
 $\log\Delta(z;J(n;\sigma_\epsilon))$ at infinity
 converge to the corresponding coefficients
 of $\log\Delta(z;\sigma_\epsilon)$ we get 
$$
 H_A(J_\epsilon(n))\to \Lambda_A(J_\epsilon),
\quad n\to\infty.
 $$
 Evidently $\Lambda_A(J_\epsilon)\le \Lambda_A(J)$. Therefore
 for every $\delta$ there exists $n_0$ such that 
 $$
 H_A(J_\epsilon(n))\le \Lambda_A(J)+\delta
 $$
 for all $n\ge n_0$. Since in the case under consideration 
$H_A$ is (basically) a series
 with positive terms, we get that every partial sum
 is bounded
 $$
 H_A^N(J_\epsilon(n))\le \Lambda_A(J)+\delta.
 $$
 Note that  the left--hand side does not depend
 on $n$ if $n$ is big enough. Thus
 $$
 H_A^N(J_\epsilon)\le \Lambda_A(J).
 $$
 Now, for a fixed $N$ let us pass to the limit as 
 $\epsilon\to 0$. Due to \eqref{p2}
 and continuity of $g_A$, for all $N$
 $$
 H_A^N(J)\le \Lambda_A(J).
 $$
 But this means that
 $$
 \limsup H_A(J(n))=\limsup \Lambda_A(J(n))\le \Lambda_A(J).
 $$
 Using Corollary \ref{c2.3} we get
 $$
 H_A(J)=\lim H_A(J(n))=\lim \Lambda_A(J(n))= \Lambda_A(J).
 $$

Finally,  starting with the condition that  series 
\eqref{s1} converges we conclude that
 $\limsup H_A(J(n))=\limsup \Lambda_A(J(n))<\infty$. Therefore,
 due to Corollary \ref{c2.3}, we have  $\Lambda_A(J)<\infty$
 and this completes the proof.
 
 \section{Asymptotic of orthonormal polynomials}

\begin{proof}[Proof of Theorem \ref{A}]
First let us mention that simultaneously with
the convergence
$$
\Lambda(J(n))=\int d\lambda_n\to
\Lambda(J)=\int d\lambda,
$$
we proved 
\begin{equation}\label{last1}
\lim_{n\to\infty}\int P(x) d\lambda_n(x)
=\int P(x) d\lambda(x)
\end{equation}
for every $P(x)=Q^2(x)$ and hence \eqref{last1}
holds for all polynomials.
Since the variations of $\lambda_n$'s are uniformly bounded
and since there is a finite interval $[\alpha_1,\alpha_2]$ containing
the support of each measure $\lambda_n$ in the family, $\lambda_n$
converges weakly to
$\lambda$.

We will  estimate the difference
$$
\left|
\int\frac{d\lambda_n}{x-z}-\int\frac{d\lambda}{x-z}
\right|
$$
on a system of contours of the form
$$
\tau= \{z=x+iy: a\le x\le b,\  y=\pm c;\quad 
|y|\le c, \ x=a,b\}
$$
that shrink to the interval $[-2,2]$.

Integrating by parts, on a horizontal line we have
\begin{equation*}
\begin{split}
\left|
\int\frac{(\lambda-\lambda_n)\,dx}{(x-z)^2}
\right|
\le&
\frac{\int_{\alpha_1}^{\alpha_2}|\lambda-\lambda_n|\,dx}
{c^2}+
|\lambda(\alpha_2)-\lambda_n(\alpha_2)|
\int_{\alpha_2}^{\infty}
\frac{dx}{|x-z|^2}\\
\le&
\frac{\int_{\alpha_1}^{\alpha_2}|\lambda-\lambda_n|\,dx}
{c^2}+
\frac{|\lambda(\alpha_2)-\lambda_n(\alpha_2)|}{c}.
\end{split}
\end{equation*}
Since the $\lambda_n(x)$ are uniformly bounded and
$\lim_{n\to\infty}\lambda_n(x)=\lambda(x)$ for all $x$, the above
estimate shows that
for every $\epsilon>0$ there exists $n_0$ such that
$$
\left|
\int\frac{d\lambda_n}{x-z}-\int\frac{d\lambda}{x-z}
\right|\le \epsilon,\ n\ge n_0,
$$
when $z$ runs on a horizontal line of the contour $\tau$.

Next, let us consider, say, the right vertical line on $\tau$.
Assume that $b$ is between of two consequent points $x_{k+1}<x_k$ of 
the set $X$. We can even specify $b=(x_{k+1}+x_k)/2$. The point
is that starting with a suitable $n$ the interval
$[b-\delta/2,b+\delta/2]$ is in a gap of the support of 
$\lambda-\lambda_n$. Here $\delta:=(x_k-x_{k+1})/2$.
Put $\tilde\lambda(x)=\lambda(x)-\lambda(b)$
and $\tilde\lambda_n(x)=\lambda_n(x)-\lambda_n(b)$.
Doing basically the same as on a horizontal line, we get
\begin{equation*}
\begin{split}
\left|
\int_{b+\delta/2}^{\infty}\frac{(\tilde\lambda-
\tilde\lambda_n)\,dx}{(x-z)^2}
\right|
\le&
\frac{\int_{b+\delta/2}^{\alpha_2}|\tilde\lambda-
\tilde\lambda_n|\,dx}
{(\delta/2)^2}+
|\tilde\lambda(\alpha_2)-\tilde\lambda_n(\alpha_2)|
\int_{\alpha_2}^{\infty}
\frac{dx}{|x-z|^2}\\
\le&
\frac{\int_{\alpha_1}^{\alpha_2}|\tilde\lambda-
\tilde\lambda_n|\,dx}
{(\delta/2)^2}+
\frac{|\tilde\lambda(\alpha_2)-\tilde\lambda_n(\alpha_2)|}{(\delta/2)},
\end{split}
\end{equation*}
and the same estimation for  $\int_{-\infty}^{b-\delta/2}$.

In other words the estimation
\begin{equation}\label{last5}
\left|A(z)\sqrt{z^2-4}\log\Delta_n(z) -B_n(z)
-\int\frac{d\lambda}{x-z}
\right|\le {\epsilon}
\end{equation}
holds on  the rectangle $\tau$
if $n\ge n_0$.

Introduce the holomorphic function
$D(z)$ by \eqref{D}, 
$z\in\overline \CC\setminus [-2,2]$,
and consider the difference
$$
\left|\Delta_n(z)e^{-\frac{B_n(z)}{A(z)\sqrt{z^2-4}}}-D(z)
\right|=
|D(z)|
\left|
e^
\frac{A(z)\sqrt{z^2-4}\log\Delta_n(z) -B_n(z)
-\int{\frac{d\lambda}{x-z}}}{A(z)\sqrt{z^2-4}}
-1
\right|
$$
on the contour $\tau$. Due to \eqref{last5} the difference is
uniformly small on the contour and therefore also in the exterior
of the rectangle.

Thus we have
\begin{equation}\label{last2}
\zeta^n (p_n P_n(z)-\zeta P_{n-1}(z))
\exp\left(-\frac{ B_n(z)}{A(z)\sqrt{z^2-4}}
\right)\to D(z)
 \end{equation}
 uniformly in  the domain $\overline \CC\setminus
 [-2,2]$.
 Let us derive from this an asymptotic for the orthonormal
 polynomials properly.
 
 First of all due to \eqref{p1} 
we have \cite{R}
$$
 \frac{P_{n-1}(z)}{p_n P_{n}(z)}\to\zeta
 $$
 uniformly in $\overline \CC\setminus
 [-2,2]$.
 Therefore from \eqref{last2} we get
 \begin{equation}\label{last3}
\zeta^n P_n(z)
\exp\left(-\frac{ B_n(z)}{A(z)\sqrt{z^2-4}}
\right)\to \frac{D(z)}{1-\zeta^2}.
 \end{equation}
 Next we will  adjust a bit the polynomials $B_n$ in \eqref{last3}.

 Let $\tilde J(n)$ be $n\times n$ matrix with coefficients
 $p_k$, $q_k$, respectively $\tilde J_0(n)$ is $n$ by $n$
 matrix that we obtain cutting the Chebyshev matrix $J_0$.
Recall that
 $$
 P_n(z)=\frac{1}{p_1...p_n}\det(z-\tilde J(n))
 $$
 in particular
 $$
\det(z-\tilde J_0(n))=\frac{\zeta^{-n-1}-\zeta^{n+1}}
{\zeta^{-1}-\zeta}.
 $$
 That is
 $$
 \frac{1}{p_1...p_n}\frac{\det(z-\tilde J(n))}
 {\det(z-\tilde J_0(n))}=
 (\zeta^{-1}-\zeta)\frac{\zeta^{n+1} P_n(z)}{1-\zeta^{2n+2}},
 $$
 and hence
 \begin{equation*}
 \begin{split}
  &\log(\zeta^{n+1}\sqrt{z^2-4} P_n(z))\\&=
  -\log(p_1... p_n)-\frac{\tr(\tilde J(n)-\tilde J_0(n))} z
  -\frac {\tr(\tilde J^2(n)-\tilde J_0^2(n))}{2z^2}-...\ .
  \end{split}
\end{equation*}
 
 Thus we can substitute $B_n(z)$
 by the polynomial $\tilde B_n(z)$, which is
 uniquely defined by 
 $$
 \log(\zeta^{n+1}\sqrt{z^2-4} P_n(z))
 -\frac{\tilde B_n(z)}{A(z)\sqrt{z^2-4}}
=\underline O\left(\frac 1{z^{m+2}}\right),
$$
 since by condition \eqref{p1}
 for any fixed $k$
 $$
 \tr(J^k(n)-J_0^k)-
  \tr(\tilde J^k(n)-\tilde J_0^k(n))\to 0,
  \ n\to\infty.
  $$
 \end{proof}

\section {Appendix: Laptev--Naboko--Safronov Example}

It is more convenient (uniform) to use two sided Jacobi
matrices acting in
$l^2(\ZZ)$. In particular,  then the function $H_A(J)$
is positive. 

\subsection{Positive definite Hankel minus Toeplitz}
Recall that the
Chebyshev polynomials of the second
kind
$U_l(z)$  form an orthogonal system with respect to
the weight
$\sqrt{4-x^2}$,
\begin{equation}\label{aaa1}
\frac 1 \pi\int_{-2}^2 U_l(x) U_k(x)\sqrt{4-x^2}\,dx=2\delta_{k,l},
\end{equation}
where
\begin{equation}\label{aaa2}
U_l(z):=\frac{\zeta^{-l}-\zeta^{l}}{\zeta^{-1}-\zeta}, \quad
z=\zeta^{-1}+\zeta.
\end{equation}
Note also that the following map transforms
the polynomials of the second kind into
the
Chebyshev polynomials of the first
kind
\begin{equation}\label{aaa3}
z U_l(z) -\frac 1 \pi\int_{-2}^2
\frac{U_l(x)-U_l(z)}{x-z}\sqrt{4-x^2}\,dx
=T_l(z).
\end{equation}
\begin{lemma} \label{aaat1}
For $m\not=n$
\begin{equation}\label{aaa4}
H_{U_m U_n}(J)=
\tr\left\{\frac{T_{m+n}}{m+n}-
\frac{T_{|m-n|}}{|m-n|}
\right\}^J_{J_0},
\end{equation}
and
\begin{equation}\label{aaa5}
H_{U_n^2}(J)=
\tr\left\{\frac{T_{2n}}{2n}-
\sum_i \log p_i^2
\right\}^J_{J_0}=
\tr\left\{\frac{T_{n}^2}{2n}-
\sum_i \log p_i^2
\right\}^J_{J_0}.
\end{equation}
\end{lemma}
\noindent

\begin{proof}
We have
\begin{equation*}
\begin{split}
\Phi'(z)=zU_m(z) U_n(z)-&
\frac 1\pi\int U_m(x)\frac{U_n(x)-U_n(z)}{x-z}
\sqrt{4-x^2}dx\\
-&
\frac 1\pi\int \frac{U_m(x)-U_m(z)}{x-z}U_n(z)
\sqrt{4-x^2}dx.
\end{split}
\end{equation*}
Using \eqref{aaa1}, \eqref{aaa2}, \eqref{aaa3} we have for $m>n$
\begin{equation*}
\begin{split}
\Phi'(z)=&zU_m(z) U_n(z)-
\frac 1\pi\int \frac{U_m(x)-U_m(z)}{x-z}
\sqrt{4-x^2}dx\,U_n(z)\\
=&T_m(z) U_n(z)= U_{m+n}(z)-U_{m-n}(z).
\end{split}
\end{equation*}
Since $T_k'=kU_k$, $k\ge 1$, we get 
$$
\Phi(z)=
\frac{T_{m+n}(z)}{m+n}-
\frac{T_{m-n}(z)}{m-n}+
\rm{const}.
$$
By orthogonality also
$$
a=\frac 1\pi\int_{-2}^2 U_m(x) U_n(x)\sqrt{4-x^2}\,dx=0.
$$
Thus \eqref{aaa4} is proved. A proof of \eqref{aaa5} requires just a
minor modification.
\end{proof}

\begin{proposition}\label{aaat2}
Let $J$ be a finite dimensional perturbation of
$J_0$. Define
\begin{equation}
a_k(J)=\begin{cases}\tr\left\{\frac{T_{k}}{k}
\right\}^J_{J_0},
 &k\ge 1\\
 \sum_i \log p_i^2
& k=0
\end{cases}
\end{equation}
Then the matrix $\{a_{k+l}(J)-a_{|k-l|}(J)\}_{k\ge 1, l\ge 1}$
is positive.
\end{proposition}
\begin{proof}
Put $A=|B|^2$ with $B=\sum_l U_l c_l$. Since $H_A(J)\ge 0$,
due to Lemma \ref{aaat1}, we get
$$
\sum_{k\ge 1, l\ge 1}\{a_{k+l}(J)-a_{|k-l|}(J)\}c_k\overline{c_l}\ge 0.
$$
\end{proof}

\noindent
Note that continuous positive kernels of this kind are a classical
object, see e.g. \cite{Akh}.

\subsection{Laptev--Naboko--Safronov example: $A=U_l^2$} This case was
considered in \cite{LNS}. 
\begin{proposition}\label{aaat3}
Let $A(z)=U_l^2(z)$.
Then $\Lambda_A(J)<\infty$ if and only if
$T_l(J)-T_l(J_0)$ is Hilbert--Schmidt.
\end{proposition}

\begin{proof} Due to Lemma \ref{aaat1}
$$
H_A(J)=\tr\frac{T^2_l(J)-T^2_l(J_0)}{2l}-
2\sum\log p_i.
$$
Note that a row in the matrix $T_l(J)$ is of the 
form
\begin{equation*}
\langle e_i| T_l(J)=
\begin{bmatrix}
\hdots &0 &(t_l)_{i-l}& (\tilde q_l)_i
& (t_l)_{i}& 0& \hdots
\end{bmatrix},
\end{equation*}
where
$(t_l)_{i}=p_{i+1}p_{i+2}...p_{i+l}$
and $(\tilde q_l)_i$ is a row--vector of dimension $2l-1$.
Therefore
$$
H_A(J)=
\frac 1 l\left\{\sum\frac{(\tilde q_l)_i(\tilde q_l)^*_i}{2}+
\sum((t_l)_i^2-1-\log(t_l)_i^2)
\right\}
$$
and the condition $H_A(J)<\infty$ is equivalent to
$T_l(J)-T_l(J_0)$ is a Hilbert--Schmidt operator.
\end{proof}

It is possible to reformulate the above condition
in terms of the coefficient sequences of $J$.

\begin{theorem}\label{aaat4}
Let $A(z)=U_n^2(z)$.
Then $\Lambda_A(J)<\infty$ if and only if
\begin{equation}\label{aaa7}
\{\sum_{k=1}^n u_{j+k}\}\in l^2,\ 
\{\sum_{k=1}^n q_{j+k}\}\in l^2, \ 
\{u_j^2\}\in l^2,\ \{q_j^2\}\in l^2,
\end{equation}
where $u_j=p_j^2-1$.
\end{theorem}
A proof is splitted in several lemmas.

\begin{lemma}\label{lns1}
Let $J=S^{-1}\PP+\QQ+\PP S$ and
$$
T_n(J)=\{...+\Lambda_0(n)+\Lambda_1(n)S+...+\Lambda_n(n)S^n\},
$$
where $\QQ,\PP,\Lambda_k(n)$ are diagonal matrices.
Then
\begin{align}
\Lambda_n(n)=&\PP\PP^{(-1)}...\PP^{(-n+1)} \label{ls1}\\
\Lambda_{n-1}(n)=&
\PP...\PP^{(-n+2)}\{\QQ+\QQ^{(-1)}+...+\QQ^{(-n+1)}\}
\label{ls2}
\end{align}
and
\begin{equation}\label{ls3}
\begin{split}
\Lambda_{n-2}(n)=&
\PP...\PP^{(-n+3)}
\{[(\PP^{(1)})^2-I+\PP^2-I+
...+(\PP^{(-n+3)})^2-I]\\
&\begin{matrix}
 +\QQ[\QQ+&
\QQ^{(-1)}+&
...&+\QQ^{(-n+2)}]\\
   &+\QQ^{(-1)}[\QQ^{(-1)}+&...&+\QQ^{(-n+2)}]\\
   &  & +&...\quad \\
   & & &+\QQ^{(-n+2)}\QQ^{(-n+2)}\}.
\end{matrix}
\end{split}
\end{equation}
\end{lemma}
\begin{proof}
All three formulas can be proved by induction using
$$
T_n(J)=JT_{n-1}(J)-T_{n-2}.
$$
Let us prove \eqref{ls3}. We have
\begin{equation*}
\begin{split}
\Lambda_{n-2}(n)=
S^{-1}\PP \Lambda_{n-1}(n-1)S
+&\QQ \Lambda_{n-2}(n-1)
+\PP S \Lambda_{n-3}(n-1)S^{-1}\\
-&\Lambda_{n-2}(n-2).
\end{split}
\end{equation*}
Substituting \eqref{ls1} and \eqref{ls2} we get
\begin{equation*}
\begin{split}
\Lambda_{n-2}(n)=&
S^{-1}\PP \PP\PP^{(-1)}...\PP^{(-n+2)}S\\
+&\QQ \PP...\PP^{(-n+3)}\{\QQ+\QQ^{(-1)}+...+\QQ^{(-n+2)}\}\\
+&\PP S \Lambda_{n-3}(n-1)S^{-1}
-\PP\PP^{(-1)}...\PP^{(-n+3)}\\
=&\PP\PP^{(-1)}...\PP^{(-n+3)}
\{(\PP^{(1)})^2-I+
\QQ[\QQ+\QQ^{(-1)}+...+\QQ^{(-n+2)}]\\
+&\PP\Lambda^{(-1)}_{n-3}(n-1)\}.
\end{split}
\end{equation*}
Iterating the last relation we obtain \eqref{ls3}.
\end{proof}

\begin{lemma}
If $T_n(J)-T_n(J_0)$ is Hilbert--Schmidt then
relations \eqref{aaa7} are fulfilled.
\end{lemma}

\begin{proof}
Since $\Lambda_{n}(n)-I$,
$\Lambda_{n-1}(n)$ and $\Lambda_{n-2}(n)$
are Hilbert--Schmidt operators, using Lemma \ref{lns1}, we have 
\begin{equation}\label{ls4}
\{p_{1+i}...p_{n+i}-1\}\in l^2
\end{equation}
\begin{equation}\label{ls5}
\{p_{1+i}...p_{n-1+i}(q_{i}+...+q_{n-1+i})\}\in l^2
\end{equation}
and
\begin{equation}\label{ls6}
\left\{p_{1+i}...p_{n-2+i}\left[\sum^{i+n-1}_{k=i}(p^2_{k}-1)
+\sum^{i+n-2}_{k=i}q_k^2
+\sum_{i\le k<l\le i+n-2}q_k q_l\right]\right\}
\in l^2.
\end{equation}
Having in mind \eqref{ls4}  we simplify
\eqref{ls5} and \eqref{ls6}
$$
\{q_{i}+...+q_{n-1+i}\}\in l^2
$$
and
\begin{equation}\label{ls7}
\left\{\sum^{i+n-1}_{k=i}(p^2_{k}-1)
+\frac 1 2\sum^{i+n-2}_{k=i}q_k^2
+\frac 1 2
\left(\sum_{k=i}^{i+n-2}q_k\right)^2\right\}
\in l^2.
\end{equation}
Now we wish to separate ``$p$" and ``$q$" conditions in
\eqref{ls7}. It is evident that $a+b\in l^2$
implies $a\in l^2$ and $b\in l^2$ if only $a_i\ge 0$
and $b_i\ge 0$. Note that
\eqref{ls4} implies
$\{(p_{1+i}...p_{n+i})^{2/n}-1\}\in l^2$.
Thus using this condition and the inequality
$$
\frac{p_{1+i}^2+...+p_{n+i}^2-n}{n}\ge
(p_{1+i}...p_{n+i})^{2/n}-1
$$
we get from \eqref{ls7} $\{q_i^2\}\in l^2$ and
$\{\sum^{n}_{k=1}(p_{i+k}^2-1)\}\in l^2$.

Finally we note that
\begin{equation*}
\begin{split}
(p_1-1)^2+...+(p_n-1)^2=&
(p_1^2-1)+...+(p^2_n-1)\\
-&2\{(p_1-1)+...+(p_n-1)\}.
\end{split}
\end{equation*}
Since
\begin{equation*}
\begin{split}
2n\{(p_1...p_n)^{1/n}-1\}
\le&2\{(p_1-1)+...+(p_n-1)\}\\
\le&
(p_1^2-1)+...+(p^2_n-1)
\end{split}
\end{equation*}
we have $\{\sum_{k=1}^n (p_{i+k}-1)\}\in l^2$ and therefore
$\{(p_{i}-1)^2\}\in l^2$.

\end{proof}

The following lemma can be shown by induction.
\begin{lemma}\label{aaat7}
Let $J=J_0+dJ$ then
\begin{equation}
dT_l(J)e_0=
\sum_{k=0}^{l-1}S^{1-l}S^k[dJ+...+dJ^{(1-l)}]S^ke_0=
\begin{bmatrix}
0\\
d p_{-l+1}+...+d p_{0}\\
d q_{-l+1}+...+d q_{0}\\
2d p_{-l+2}+...+2d p_{1}\\
d q_{-l+2}+...+d q_{1}\\
2d p_{-l+3}+...+2d p_{2}\\
\vdots\\
d p_{1}+...+d p_{l}\\
0
\end{bmatrix}.
\end{equation}
\end{lemma}

\begin{proof}[Proof of the Theorem \ref{aaat4}]
We only have to show that conditions \eqref{aaa7} imply $T(J)-T(J_0)$
is Hilbert--Schmidt. Note that each entry is a polynomial
of $q_j, u_i$ with $u_i=p_i-1$. Moreover, the linear term is described
in Lemma \ref{aaat7}. Note also that  the sequences 
$\{u_i^l q_{i+j}^k\}_i$,
$\{u_i^l u_{i+j}^k\}_i$,
$\{q_i^l q_{i+j}^k\}_i$ 
belong to $l^2$ for $k+l\ge 2$. Thus, having in
mind the structure of the matrix $T(J)-T(J_0)$, we get that each
diagonal forms an $l^2$--sequence, as was to be proved.
\end{proof}

\subsection{Simon's conjecture}
Since  $H_A(J_0)=0$ 
and $H_A(J)\ge 0$
the decomposition of $H_A$ about $J_0$ begins with a
quadratic form, more exactly:
\begin{lemma}
Let $J=J_0+dJ$ then
the decomposition of $H_A$ about $J_0$ begins with 
\begin{equation}\label{l041}
H_A(J)=\frac 1 2\langle dj| A(J_0)|dj\rangle+\dots
\end{equation}
where  $\langle dj|=\{\dots, 2d p_0, d q_0, 2d p_1, dq_1,\dots\}$. 
\end{lemma}

\begin{proof}
We start with the formula
$$
d H_A(J)=\tr \{A(J){\rm Re}(Z^{-1}-Z)\,dJ\},
$$
where $Z$ is the lower triangle solution of the equation
$Z^{-1}+Z=J$.
Note that the decomposition of $Z^{-1}-Z$
 about $J_0$ is of the form
$$
Z^{-1}-Z=S^{-1}-S+dJ-2dZ+...\ .
$$
Using
$$
dJ=-Z^{-1}dZ Z^{-1}+ dZ
$$
we get
$$
-dZ|_{Z=S}=
[Z dJ Z + Z(-dZ)Z]_{Z=S}=
SdJS+S^2dJS^2+...\ .
$$
Therefore the leading term in the decomposition
of  Re$(Z^{-1}-Z)$ is the Hankel operator
$$
\Gamma=...+S^{-1}dJS^{-1}+dJ+SdJS+... ,
$$
and
$$
H_A(J)=\frac 1 2 \tr\{A(J_0)\Gamma\,dJ\}+...\ .
$$

Let us mention that $\Gamma e_0=dj$, thus we can rewrite
this Hankel operator into the form
$$
\Gamma=\sum S^k|dj\rangle\langle e_0| S^k.
$$
Since $A(J_0)$ and $S$ commute and $\Gamma S=S^{-1}\Gamma$ we get
\begin{equation*}
\begin{split}
\tr\{A(J_0)\Gamma\,dJ\}&=\tr\{A(J_0)\Gamma(S^{-1}d\PP+d\QQ+d\PP S)\}\\
&=\tr\{A(J_0)\Gamma(2S^{-1}d\PP+d\QQ)\}.
\end{split}
\end{equation*}
Substituting $\Gamma$ we obtain
\begin{equation*}
\begin{split}
\tr\{A(J_0)\Gamma\,dJ\}&=\tr\{A(J_0)
(\sum S^k|dj\rangle\langle e_0| S^k)
(2S^{-1}d\PP+d\QQ)\}\\
&=\tr\{A(J_0) |dj\rangle\langle e_0| \sum (2S^{k-1}
d\PP S^k+S^{k}d\QQ S^k)\}.
\end{split}
\end{equation*}
But
$\langle e_0| \sum (2S^{k-1}
d\PP S^k+S^{k}d\QQ S^k)=\langle dj|$
and this completes the proof.
\end{proof}

We believe that related to
this quadratic form condition
\begin{equation}\label{r2}
\langle A(J_0) dj,dj\rangle<\infty,
\end{equation}
 should play an important
role in a counterpart of Simon's conjecture formulated for the unit
circle  in several talks, for
example \cite{S}. 
Specifically, in Laptev--Naboko--Safronov case,
where
$$
A(J_0)=(I+S^2+...+S^{2l-2})^*(I+S^2+...+S^{2l-2}),
$$ 
condition \eqref{r2} means
\begin{equation*}
\begin{split}
&\{d q_{i+1}+d q_{i+2}+...+d q_{i+l}\}\in l^2(\ZZ),\\
&\{2d p_{i+1}+2d p_{i+2}+...+
2d p_{i+l}\}\in l^2(\ZZ),
\end{split}
\end{equation*}
compare \eqref{aaa7}.


\bibliographystyle{amsplain}

\end{document}